\documentclass[11pt,letterpaper,reqno]{amsart}
\usepackage[utf8]{inputenc}
\usepackage{dsfont,amsmath,amsfonts,amssymb,amsthm,tikz,graphicx,enumitem,polynom,mathrsfs,mathabx,mathtools, verbatim}

\usetikzlibrary{shapes.geometric}
\usepackage[left=1in, right=1in, bottom=1in, top=1in, bottom=1in]{geometry}

\newcommand{\z}{\mathbb{Z}}
\newcommand{\rl}{\mathbb{R}}
\newcommand{\q}{\mathbb{Q}}
\newcommand{\n}{\mathbb{N}}

\DeclareMathOperator{\slz}{\mathrm{SL}_2(\z)}
\DeclareMathOperator{\slt}{\mathrm{SL}_2}
\newcommand{\gl}{\textrm{GL}_2}
\newcommand{\mkh}{\widehat{M}_k^{(\infty)}}
\newcommand{\mkhd}{\widehat{M}_{2-k}^{(\infty)}}
\newcommand{\mki}{M_k^{(\infty)}}

\DeclareMathOperator{\tr}{\mathrm{Tr}}
\newtheorem{theorem}{Theorem}
\newtheorem{proposition}{Proposition}

\newtheorem*{theorem*}{Theorem}
\newtheorem{innercustomgeneric}{\customgenericname}
\providecommand{\customgenericname}{}
\newcommand{\newcustomtheorem}[2]{%
  \newenvironment{#1}[1]
  {%
   \ifdefined\crefalias\crefalias{innercustomgeneric}{#2}\fi
   \renewcommand\customgenericname{#2}%
   \renewcommand\theinnercustomgeneric{##1}%
   \innercustomgeneric
  }
  {\endinnercustomgeneric}%
  \ifdefined\crefname\crefname{#2}{#2}{#2s}\fi
}
\newcustomtheorem{customthm}{Theorem}

\mathtoolsset{showonlyrefs=true}
\numberwithin{proposition}{section}
\numberwithin{figure}{section}
\numberwithin{table}{section}
\numberwithin{theorem}{section}

\usepackage{etoolbox}
\usepackage{hycolor} 
\usepackage[colorlinks=true]{hyperref}
\hypersetup{
  colorlinks=true,
  linkcolor=blue,
  citecolor=blue,
  linkbordercolor={0 0 1}
}
\title[The Trace Operator on Modular Grids]{The Effect of the Trace Operator on the Duality of Modular Grids in Genus Zero Levels}

\author{Archer Clayton}
\address{Department of Mathematics, Brigham Young University, Provo, UT 84602}
\email{ac727@byu.edu}
\author{Paul Jenkins}
\address{Department of Mathematics, Brigham Young University, Provo, UT 84602}
\email{jenkins@math.byu.edu}
\date{\today}
\begin{document}
\begin{abstract}
    Griffin, the second author, and Molnar studied coefficient duality for canonical bases for a broad range of spaces of weakly holomorphic modular forms, showing that the Fourier coefficients of canonical basis elements appear as negatives of Fourier coefficients for elements of a canonical basis of a related space of forms. We investigate the effect of the trace operator on this duality for modular forms for $\Gamma_0(N)$ of genus zero and show exactly when duality still holds after applying the trace operator.
\end{abstract}
\maketitle

\section{Introduction}
A weakly holomorphic modular form is a modular form for which poles are allowed at cusps. Many spaces of weakly holomorphic modular forms have row-reduced canonical bases with Fourier coefficients that exhibit an intriguing duality property between weights $k$ and $2-k$. Pairs of sequences of weakly holomorphic modular forms that exhibit this duality are called \emph{modular grids}.

For example, Duke and the second author constructed canonical bases for weakly holomorphic modular forms for $\slz$ and provided an example of a level 1 modular grid \cite{duke}.
As usual, denote the space of holomorphic modular forms for $\slz$ of weight $k$ by $M_k(\slz)$, the subspace of cusp forms by $S_k(\slz)$, and the space of weakly holomorphic modular forms of weight $k$ by $M_k^!(\slz)$, and write  $q=e^{2\pi i z}$. The Fourier expansions of the first few canonical basis elements $f_i\in M_0^!(\slz)$ and canonical basis elements $g_j\in M_2^!(\slz)$ are given by
 \begin{align*}
        f_0(z)&=1+0q+0q^2+0q^3+\cdots,\\
        f_1(z)&=q^{-1}+196884q+21493760q^2+864299970q^3+\cdots,\\
        f_2(z)&=q^{-2}+42987520q+40491909396q^2+8504046600192q^3+\cdots,\\
        f_3(z)&=q^{-3}+2592899910q+12756069900288q^2+9529320689550144q^3+\cdots,\\
        &\\
        g_1(z)&=q^{-1}-196884q-42987520q^2-2592899910q^3+\cdots,\\
        g_2(z)&=q^{-2}-21493760q-40491909396q^2-12756069900288q^3+\cdots,\\
        g_3(z)&=q^{-3}-864299970q-8504046600192q^2-9529320689550144q^3+\cdots.\\
    \end{align*}
If the weight $0$ basis elements are written as
\begin{equation}
    f_m(z)=q^{-m}+\sum_{n\geq 1}a_{0}(m,n)q^n,
\end{equation}
and the weight 2 basis elements are written as
\begin{equation}
    g_n(z)=q^{-n}+\sum_{m\geq 0}a_{2}(n,m)q^m,
\end{equation}
then the Fourier coefficients satisfy
\begin{equation}
    a_{0}(m,n)=-a_{2}(n,m).
\end{equation}
In other words, the $n$th Fourier coefficient of the $m$th basis element for the weight 0 space is exactly the negative of the $m$th coefficient of the $n$th basis element for the weight 2 space.

 This phenomenon was earlier observed by Zagier \cite{Zagier} for half integer weight modular forms in the Kohnen plus space. Coefficient duality appears in many other settings as well, including harmonic Maass forms, vector valued modular forms, and related modular objects; see for example \cite{ElGuindy}, \cite{folsom}, \cite{guerzhoy}.

 Let $\Gamma$ be a congruence subgroup of $\slz$, let $k\in \frac{1}{2}\z$, and let $M_k^!(\Gamma)$ be the space of weakly holomorphic modular forms of weight $k$ for $\Gamma$. Let $\mki(\Gamma)\subseteq M_k^!(\Gamma)$ denote the space of weakly holomorphic forms for $\Gamma$ of weight $k$ that have poles only at the cusp at $\infty$, and let $\mkh(\Gamma)\subseteq\mki(\Gamma)$ be the subspace of forms that vanish at all cusps other than $\infty$. Note that \[\mki(\slz)=\mkh(\slz)=M_k^!(\slz).\]
 For notational convenience, when $\Gamma=\Gamma_0(N)$, we often write $M_k(\Gamma_0(N)),$ $M_k^!(\Gamma_0(N))$, $\mki(\Gamma_0(N)),$ and $\mkh(\Gamma_0(N))$ as $M_k(N)$, $M_k^!(N)$, $\mki(N)$, and $\mkh(N)$ respectively.

 In general, the spaces $\mki(\Gamma)$ and $\widehat{M}_{2-k}^{(\infty)}(\Gamma)$ have canonical bases that form modular grids. The existence of such canonical bases for these spaces was proven by Griffin, the second author, and Molnar in \cite{griffin}.
 \begin{theorem}\cite[Theorem~1.1]{griffin}\label{DualityTheorem}
     Let $\Gamma\subseteq \slt(\rl)$ be a subgroup that is commensurable with $\slz$ and let $k\in\frac{1}{2}\z$. There exist row-reduced canonical bases for $\mki(\Gamma)$ and $\widehat{M}_{2-k}^{(\infty)}(\Gamma)$ with elements
     \begin{equation}
         f_{k,m}(z)=q^{-m}+\sum a_{k}(m,n)q^n \quad\textrm{and}\quad g_{2-k,n}(z)=q^{-n}+\sum b_{2-k}(n,m)q^m,
     \end{equation}
     whose Fourier coefficients satisfy
     \begin{equation}
        a_{k}(m,n)= -b_{2-k}(n,m).
     \end{equation}
 \end{theorem}
\noindent In addition, they showed that these modular grids have explicit generating functions encoding them.
 \begin{theorem}\cite[Theorem~1.2]{griffin}\label{GenFunTheorem}
     Let $\Gamma, k$ be as above and let $\left(f_{k,m}(z)\right)_m$ and $\left(g_{2-k,n}(z)\right)_n$ be row-reduced canonical bases for the spaces $\mki(\Gamma)$ and $\mkhd(\Gamma)$. There is an explicit generating function $\mathcal{H}_k(z,\tau)$ which encodes the modular grid; that is, on appropriate regions of $\mathcal{H}\times\mathcal{H}$, the function $\mathcal{H}_k(z,\tau)$ satisfies
     \begin{equation}
         \mathcal{H}_k(z,\tau)=\sum_m f_{k,m}(z)e^{2\pi i m \tau}\quad\textrm{and}\quad  \mathcal{H}_k(z,\tau)=-\sum_n g_{2-k,n}(\tau)e^{2\pi i n z}.
     \end{equation}
 \end{theorem}
\noindent Note that the theorems from \cite{griffin} are more general than those given here.

Many such formulas for these generating functions for specific $\Gamma$ and $k$ were previously known.  For example, if $\Gamma = \Gamma_0(4)$ and $k=2\ell$, then \cite{haddock} gives the generating function
\[\mathcal{H}_k(z, \tau) = \frac{f_{k, -\ell}(z) g_{2-k, \ell + 1}(\tau)}{f_{0, 1}(\tau) - f_{0, 1}(z)}\]
in terms of the canonical basis elements from Theorem~\ref{DualityTheorem}.  Other similar generating functions appear in~\cite{ElGuindy}.

  In \cite[Section~7]{griffin}, new modular grids are constructed using linear operators on existing modular grids. In fact, sometimes different constructions result in the same grid. One interesting example is the effect of the Hecke operators in weights $k$ and $2-k$.
\begin{theorem}\cite[Theorem~1.3]{griffin}
Let $k$ be a positive even integer, and let $\Gamma\supseteq\Gamma_0(N)$ for some $N\in\n$ such that $S_k(\Gamma)={0}$. If $M$ is a positive integer coprime to $N$, then applying the Hecke operator $T_M$ in weight $k$ to the weight $k$ basis elements and in weight $2-k$ to the weight $2-k$ basis elements gives
\begin{equation}
    \mathcal{H}_{2-k}(z,\tau)|_{k,\tau}T_M=\mathcal{H}_{2-k}(z,\tau)|_{2-k,z}T_M.
\end{equation}
\end{theorem}
\noindent For further information about Hecke operators $T_M$, see \cite[Chapter 5]{diamond2005first}. When the operators in weights $k$ and $2-k$ give the same result as in the case above, we say that they \emph{preserve the duality of the modular grid}.

 It is natural to investigate whether other linear operators preserve the duality of modular grids. In this paper, we study the action of the trace operator $\tr^N_M$, which maps a weakly holomorphic modular form for $\Gamma_0(N)$ to a form for $\Gamma_0(M)$, where $M\mid N$. The trace operator, like the Hecke operators, can preserve the duality of modular grids. For modular forms of integer weight for $\Gamma_0(N)$ of genus zero, we give a complete characterization of the weights in which the trace preserves duality.
\begin{theorem}\label{MainTheorem}
Let $N\in\{2, 3, 4, 5, 6, 7, 8, 9, 10, 12, 13, 16, 18, 25\}$, so that $\Gamma_0(N)$ has genus zero. Let $k$ be an even integer and let $M\mid N$. Then the trace operator $\tr^N_M$ preserves duality between grids of basis elements for $\mki(N)$ and $\mkhd(N)$ for $k=0$. Additionally, if $N\in\{2, 3, 4\}$ then $\tr^N_M$ preserves duality for $k=-2$, and $\tr^2_1$ also preserves duality for $k=-4$. If $M=N,$ then $\tr^N_M$ is the identity map and preserves duality. In all other weights and levels $\tr^N_M$ does not preserve duality.
\end{theorem}
 \noindent Duality of the modular grid is only preserved by the trace operator when related spaces of holomorphic modular forms are trivial (see Proposition \ref{Prop:LvN Trace}), similar to \cite[Theorem~1.3]{griffin}.

As an example, the generating function for forms of weight $0$ for $\Gamma_0(4)$ above is \[\mathcal{H}_0(z, \tau) = \frac{g_{2, 1}(\tau)}{f_{0, 1}(\tau) - f_{0, 1}(z)}.\]  The Fourier coefficients of the first few basis elements in weights $0$ and $2$ are
\begin{align*}
    f_{0, 0}(z) &= 1, \\
    f_{0, 1}(z) &= q^{-1} + 20 q  + 0 q^2 - 62 q^3  + 0 q^4 + 216 q^5 + \cdots, \\
    f_{0, 2}(z) &= q^{-2}  + 0 q + 276 q^2  + 0 q^3 - 2048 q^4  + 0 q^5 + 11202 q^6 + \cdots, \\
    f_{0, 3}(z) &= q^{-3}  - 186 q  + 0 q^2 + 4928 q^3  + 0 q^4 - 51831 q^5 + \cdots, \\
    &\\
    g_{2, 1}(z) &= q^{-1} - 20q + 186q^3 - 1080q^5 + \cdots, \\
    g_{2, 2}(z) &= q^{-2} - 276q^2 + 4096q^4 - 33606q^6 + \cdots, \\
    g_{2, 3}(z) &= q^{-3} +62q - 4928q^3 + 86385q^5 + \cdots.
\end{align*}
Applying the weight $0$ trace operator $\tr^4_1$ to the $f_{0, n}(z)$ gives
\begin{align*}
    \tr^4_1 (f_{0, 0}(z)) &= 1, \\
    \tr^4_1 (f_{0, 1}(z)) &= q^{-1} + 196884q + 21493760q^2 + 864299970q^3 + \cdots, \\
    \tr^4_1 (f_{0, 2}(z)) &= q^{-2} + 42987520q + 40491909396q^2 + 8504046600192q^3 + \cdots, \\
    \tr^4_1 (f_{0, 3}(z)) &= q^{-3} + 2592899910q + 12756069900288q^2 + 9529320689550144q^3 + \cdots,
\end{align*}
while applying the weight $2$ trace operator $\tr^4_1$ to the $g_{2, m}(z)$ gives
\begin{align*}
    \tr^4_1 (g_{1, 1}(z)) &= q^{-1} - 196884q - 42987520q^2 - 2592899910q^3 + \cdots, \\
    \tr^4_1 (g_{2, 2}(z)) &= q^{-2} - 21493760q - 40491909396q^2 - 12756069900288q^3 + \cdots, \\
    \tr^4_1 (g_{3, 3}(z)) &= q^{-3} - 864299970q - 8504046600192q^2 - 9529320689550144q^3 + \cdots.
\end{align*}
Duality is preserved in this case; in fact, applying the trace operator to the generating function for $\Gamma_0(4)$ in either variable gives the generating function $\mathcal{H}_0(z, \tau)$ for $\slz$.

The remainder of the paper is dedicated to proving Theorem \ref{MainTheorem} and explicitly computing the action of the trace operator on generating functions for spaces of weakly holomorphic modular forms for $\Gamma_0(N)$ of genus zero. In Section 2 we give necessary definitions and explain how to construct canonical bases for the spaces $\mki(N)$ and $\mkh(N)$ for all such values of $N$ and $k\in 2\z$. In Section 3 we introduce the trace operator and explicitly calculate traces of modular forms in $\mki(N)$ and $\mkh(N)$. This allows us to prove Theorem \ref{MainTheorem} in Section 4 and explain how the trace acts on the generating functions for modular grids for $\mki(N)$ and $\mkhd(N)$ in Section 5.

\section{Construction of Canonical Bases}

Let $f$ be a modular form of weight $k$ for $\Gamma$ and let $\gamma=\begin{psmallmatrix}
    a&b\\
    c&d
\end{psmallmatrix}\in\gl^+(\q)$. The weight $k$ slash operator $|_k$ is defined \cite[p.~165]{diamond2005first} by
\begin{equation}
    f(z)|_k\gamma=(cz+d)^{-k}(ad-bc)^{k-1}f\left(\frac{az+b}{cz+d}\right).
\end{equation}

We now describe Fourier expansions of weakly holomorphic modular forms for $\Gamma_0(N)$ at the cusps of $\Gamma_0(N)$. We do so to understand the \textit{order of vanishing} of a modular form at a cusp. This is defined to be the order of the zero at the cusp, or if there is a pole at the cusp, the negative of the order of the pole at the cusp.
Following \cite[Section 2.3]{griffin}, let $f\in M_k^!(N)$. Then for $\alpha\in\slz$, the form $f|_k\alpha=f^\alpha$ is modular for $\alpha^{-1}\Gamma_0(N)\alpha$. Hence, $f^\alpha$ has a Fourier expansion about $\infty$ given by
    \begin{equation}
        f^\alpha(z)=\sum_{n\gg -\infty}a^\alpha(n)q^n,\label{eq:othercusps}
    \end{equation}
    for some $a^\alpha(n)$. We interpret \eqref{eq:othercusps} as an expansion of $f$ about the cusp $\alpha\infty$.

      It is well known that if $\beta\in\slz$ with $\beta\infty=\alpha\infty$ under the action of $\Gamma_0(N)$, then $a^\alpha(n)=\zeta_n a^\beta(n)$, where $\zeta_n$ is a root of unity.
     Thus, although the Fourier expansion of $f$ at the cusp $\alpha\infty$ is not canonical, it has a well-defined order of vanishing at $\infty$, which we define to be the order of vanishing of $f$ at the cusp $\alpha\infty$.

We combine and build on results from \cite{duke}, \cite{garthwaite}, \cite{haddock}, \cite{Iba}, \cite{thornton}, and \cite{thorntonp} to explicitly construct modular grids for $\mki(N)$ and $\mkh(N)$ for all $N$ such that $\Gamma_0(N)$ is of genus zero. Let $v_k(N)$ denote the maximal order of vanishing at $\infty$ for a form in $\mki(N)$, and note that this will be negative if there must be a pole at $\infty$. Similarly, let $u_k(N)$ be the maximal order of vanishing at $\infty$ for a form in $\mkh(N)$. When $k$ is clear from context, we write $v_k(N)=v(N)$ and $u_k(N)=u(N)$, and when both $k$ and $N$ are clear we write $v(N)=v$ and $u(N)=u$.

For $\Gamma_0(N)$ of genus zero, row-reduced canonical basis elements for $\mki(N)$ exist for each $m\geq -v$; for $\mkh(N)$, such basis elements exist for each $m\geq -u$. Denote the Fourier coefficients of these basis elements by $a_k^{(N)}(m,n)$ and $b_k^{(N)}(m,n)$ so that their Fourier expansions have the form
\begin{align} \label{eq:basisFourier}
    f_{k,m}^{(N)}(z)&=q^{-m}+\sum_{n=-v_k(N)+1}^\infty a_k^{(N)}(m,n)q^n\in\mki(N),\\
    g_{k,m}^{(N)}(z)&=q^{-m}+\sum_{n=-u_k(N)+1}^\infty b_k^{(N)}(m,n)q^n\in\mkh(N).
\end{align}
By Theorem \ref{DualityTheorem},
\begin{equation}
    a_k^{(N)}(m,n)=-b_{2-k}^{(N)}(n,m).
\end{equation}

In each $M_k^!(N)$ with $\Gamma_0(N)$ of genus zero there is a \emph{Hauptmodul}, or a meromorphic modular form of weight zero that has exactly one simple pole. We choose Hauptmoduln that are weakly holomorphic modular forms, so the pole will always be at a cusp. These can be written in terms of the Dedekind eta function, which is  $\eta(z)=q^{1/24}\prod_{n=1}^\infty(1-q^n)$.
The $f_{k,m}^{(N)}$ and the $g_{k,m}^{(N)}$ are constructed recursively in terms of a Hauptmodul with a simple pole at $\infty$ and the previous basis elements $f_{k,r}^{(N)}$ with $r<m$ to obtain a Fourier expansion of the form in \eqref{eq:basisFourier}.

To demonstrate this process, we explicitly construct the canonical bases for $\mki(5)$ and $\mkh(5)$, using methods explained in \cite{thorntonp}. We then describe how to construct such bases for all $\mki(N)$ and $\mkh(N)$ of genus zero and provide the modular forms needed to do so.

To begin our construction of the basis for level 5, we choose the Hauptmodul given by
\begin{equation}
    \psi^{(5)}(z)=\left(\frac{\eta(z)}{\eta(5z)}\right)^6=q^{-1}-6+9q+10q^2+\cdots\in \widehat{M_0}^{(\infty)}(5).
\end{equation}
In weight $k=0$, the first basis element is the constant function $f^{(5)}_{0,0}(z)=1$. In order to recursively obtain the next basis element, we multiply the previous basis element by the Hauptmodul $\psi^{(5)}(z)$ and row reduce using previous basis elements so that the Fourier expansion is of the correct form. Here are the first few basis elements for $M_0^{(\infty)}(5)$:
\begin{align*}
    f^{(5)}_{0,0}(z)&=1,\\
    f^{(5)}_{0,1}(z)&=\psi^{(5)}(z)+6&=&~q^{-1}+9q+10q^2+\cdots,\\
    f^{(5)}_{0,2}(z)&=\psi^{(5)}(z)f^{(5)}_{0,1}(z)+6f^{(5)}_{0,1}(z)-18&=&~q^{-2}+20q+21q^2+\cdots,\\
    f^{(5)}_{0,3}(z)&=\psi^{(5)}(z)f^{(5)}_{0,2}(z)+6f^{(5)}_{0,2}(z)-9f^{(5)}_{0,1}(z)-30&=&~q^{-3}-90q+288q^2+\cdots.
\end{align*}

Now let $k$ be any even integer. The first canonical basis element of $\mki(5)$ is the form with leading coefficient 1 and a zero of highest order, or pole of lowest order, at the cusp $\infty$. In general, we write $k=4\ell+k'$ with $k' \in \{0, 2\}$, so that $\dim M_k(5)=2\ell+1$ when $k\geq 0$. If $k$ is negative, then $\ell$ will also be negative and the first basis element will have a pole at $\infty$. To construct the first basis element, we use the modular forms
\begin{align*}
    &f^{(5)}_{2,0}(z)=\frac{5E_2(5z)-E_2(z)}{4}=1+6q+18q^2+24q^3+\cdots,\\
    &f^{(5)}_{4,-2}(z)=\frac{\eta(5z)^{10}}{\eta(z)^2}=q^2+2q^3+5q^4+10q^5+\cdots.
\end{align*}
Here $E_2(z)$ is the weight 2 Eisenstein series, with Fourier expansion at $\infty$ given by
\begin{equation}
    E_2(z)=1-24\sum_{n=1}^\infty\sigma_1(n)q^n,
\end{equation}
where $\sigma_1(n)$ is the sum of divisors function (see \cite[p.~21]{diamond2005first}).
The first basis element of $\mki(5)$ is then given by
\begin{equation}
    f^{(5)}_{k,-2\ell}(z)=f^{(5)}_{4,-2}(z)^\ell f^{(5)}_{2,0}(z)^{k'/2}= q^{2\ell}+O(q^{2\ell+1}).
\end{equation}

To recursively construct the next element of the basis in weight $k$, we multiply the previous basis element by the Hauptmodul and then row reduce using previous basis elements to create the largest gap possible between the first and second nonzero terms of the Fourier expansion. In general, the basis elements are of the form
\begin{equation}
    f^{(5)}_{k,m}(z)=q^{-m}+O(q^{2\ell+1}).
\end{equation}

Now consider $\widehat{M}_2^{{(\infty)}}(5)$. We start with the first basis element from $M_2^{{(\infty)}}(5)$, which is $f^{(5)}_{2,0}(z)$. The Hauptmodul $\psi^{(5)}(z)$ has a simple zero at the cusp 0, so $f^{(5)}_{2,0}(z)\psi^{(5)}(z)$ has a simple zero at the cusp 0 and so $f^{(5)}_{2,0}(z)\psi^{(5)}(z)\in\widehat{M}_2^{{(\infty)}}(5)$. This is the first basis element of this space and has Fourier expansion given by
\begin{equation}
    g^{(5)}_{2,1}(z)=q^{-1}-9q-20q^2+\cdots.
\end{equation}
We recursively compute additional basis elements as before; the first few are
\begin{align}
    &g^{(5)}_{2,1}(z)=q^{-1}-9q-20q^2+90q^3+\cdots,\\
    &g^{(5)}_{2,2}(z)=q^{-2}-10q-21q^2-288q^3+\cdots,\\
    &g^{(5)}_{2,3}(z)=q^{-3}+30q-192q^2-144q^3+\cdots.
\end{align}

We now explain how to construct the canonical bases for $\mki(N)$ for all $N$ where $\Gamma_0(N)$ has genus zero. The appendix gives a Hauptmodul for each $N$ in terms of $\eta(z)$. Dimension formulas for spaces of holomorphic modular forms are used to calculate $v_k(N)$ for each integer weight $k$. Once $v_k(N)$ is determined, we then find modular forms with maximal order of vanishing at $\infty$ and use these to construct the first basis element in each even weight $k$. The values of $v_k(N)$ as well as the requisite modular forms are also given in the appendix. To obtain subsequent basis elements, we multiply the first by a power of the Hauptmodul and row reduce using previous basis elements, just as in level 5.

 As elements of $\mkh(N)$ vanish at all cusps away from the cusp $\infty$, the first basis element must as well. To determine the first basis element, we find a weight zero weakly holomorphic form which has a simple zero at each of the cusps away from $\infty$ and has a pole only at $\infty$ and multiply it by the first basis element of $\mki(N)$. For each $\Gamma_0(N)$ of genus zero, a modular form in $M_0^!(N)$ that has a simple zero at each cusp away from $\infty$ can be written as a polynomial in the Hauptmodul and may be computed as in \cite{thornton} by finding the value of the Hauptmodul at each cusp. The other basis elements are obtained recursively by multiplying by the Hauptmodul and row reducing using previous basis elements. The requisite modular forms for constructing these bases and the values of $u_k(N)$ are also provided in the appendix.

\section{Calculating the Trace}
 The trace operator takes a weakly holomorphic modular form $f$ in $M_k^!(N)$ and maps it to a form in $M_k^!(M)$, where $M\mid N$. It is defined by
\begin{equation}
    \tr^N_M(f)=\sum_{\gamma\in \Gamma_0(N)\backslash\Gamma_0(M)}f|_k\gamma,
\end{equation}
where $\Gamma_0(N)\backslash\Gamma_0(M)$ is a set of coset representatives for $\Gamma_0(N)$ in $\Gamma_0(M)$.

Note that in the special case where $N=p$ is prime and $M=1$, we can write the trace using the $U_p$ and $W_p$ operators to simplify computations. The $U_p$ operator is as in \cite[p.~70-71]{kilford2015modular}, and the $W_p$ operator is the Fricke involution \cite{lmfdb:fricke}. Explicit matrix computations give
\begin{equation}
    \tr^p_1: M_k^!(p)\rightarrow M_k^!(1), \quad
    \tr^p_1(f)= f+p^{2-k}f|_kW_p|_kU_p.
\end{equation}

\subsection{Calculating Traces}

By using the bases constructed in the previous section, traces are especially straightforward to compute for modular forms of certain weights. Calculating these traces suffices to prove Theorem \ref{MainTheorem}. In general, for $f\in \mki(N)$, calculating $\tr^N_M(f)$ is dependent on the \emph{principal part} of $f$, which consists of the terms of the Fourier expansion of $f$ with negative powers of $q$. This is because the principal part of $f$ determines the principal part of $\tr^N_M(f)$, and for certain weights the principal part of a modular form of weight $k$ for $\Gamma_0(N)$ completely determines the form.
    \begin{proposition}\label{Prop:LvN Trace}
        Let $N$ be a level of genus zero and let $M$ be a level such that $M\mid N$. Let $k\in 2\z$, let $f\in\mki(N)$, and let $g\in\mkh(N)$. Then $\tr^N_M(f)$ is completely determined by the principal part of $f$ when $M_k(M)=0$, and $\tr^N_M(g)$ is completely determined by the principal part and constant term in the Fourier expansion of $g$ when $S_k(M)=0$.
    \end{proposition}
The space $M_k(M)$ is trivial when $M=1$ and $k=2$ or $k<0$, and for other $M$ when $k<0$. This means that for these weights any element of $\mki(M)$ is completely determined by the principal part of its Fourier expansion at $\infty$. The space $S_k(M)$ is trivial when $M=1$ and $k=14$ or $k<12$, for $M=2$ when $k<8$, for $M=3$ when $k<6$, and for other $M$ where $\Gamma_0(M)$ has genus zero when $k<4$. In these cases, the principal part and constant term of the Fourier expansion of an element of $\mkh(M)$ completely determine the form.
    \begin{proof}[Proof of Proposition 3.1]
        Suppose that $N$ is a level of genus zero and $M$ is a positive divisor of $N$ not equal to $N$.
    Suppose that $M_k(M)=0$, and let $f\in\mki(N)$. We have
    \begin{equation}
    \tr^N_M(f)=\sum_{\gamma\in\Gamma_0(N)\backslash\Gamma_0(M)}f|_k\gamma=f+\sum_{\substack{\gamma\in\Gamma_0(N)\backslash\Gamma_0(M)\\\gamma\neq I}}f|_k\gamma \in M_k^!(M).
    \end{equation}
    Note that if $\tr^N_M(f)\in \mki(M)$, it will be uniquely determined by its principal part since there are no holomorphic forms in $\mki(M)$.
    Since $f \in \mki(N)$, it is holomorphic at every cusp of $\Gamma_0(N)$ except possibly $\infty$, and $f|_k \gamma$ must be holomorphic at every cusp of $\Gamma_0(N)$ except possibly $\gamma\infty$.  Explicit computations show that if $\gamma\neq I$ is a coset representative for $\Gamma_0(N)$ in $\Gamma_0(M)$, then under the action of $\Gamma_0(N)$, the cusp $\gamma\infty \neq \infty$, while under the action of $\Gamma_0(M)$, the cusp $\gamma\infty=\infty$.
    Since $\tr^N_M(f)\in M_k^!(M)$ and the $f|_k \gamma$ are all holomorphic at every cusp of $\Gamma_0(N)$ that is not equivalent to $\infty$ under the action of $\Gamma_0(M)$, we conclude that $\tr^N_M(f)$ must be in the space $\mki(M)$.

   Consider $\tr^N_M(f)$ as a modular form in the larger space $M_k^!(N)$. The principal part of its Fourier expansion at $\infty$ will be the sum of the principal part of $f$ and the principal parts of the $f|_k \gamma$.  But since $f$ is holomorphic at each cusp of $\Gamma_0(N)$ except possibly $\infty$, each $f|_k\gamma$ is holomorphic at $\infty$ and has no principal part. Thus, the principal part of $f$ determines the principal part of the Fourier expansion of $\tr^N_M(f)$.

   Because $\tr^N_M(f)\in \mki(M)$ and $M_k(M) = 0$, the form $\tr^N_M(f)$ is completely determined by its principal part.  Thus, $\tr^N_M(f)$ is equal to the linear combination of basis elements from $\mki(M)$ with principal part that equals the principal part of $f$.

    Similarly, suppose that $S_k(M)=0$ and let $g\in\mkh(N)$. Then
    \begin{equation}
    \tr^N_M(g)=\sum_{\gamma\in\Gamma_0(N)\backslash\Gamma_0(M)}g|_k\gamma\in M_k^!(M).
    \end{equation}
If $\tr^N_M(g)\in\mkh(M)$, it will be uniquely determined by its principal part and constant term since there are no cusp forms in $\mkh(M)$. Since $g\in \mkh(N)$, it vanishes at every cusp of $\Gamma_0(N)$ except possibly $\infty$, and as above this means that $g|_k\gamma$ vanishes at every cusp of $\Gamma_0(N)$ except possibly $\gamma\infty$. If $\gamma\neq I$ is a coset representative for $\Gamma_0(N)$ in $\Gamma_0(M)$, it permutes the cusps of $\Gamma_0(N)$ and $\Gamma_0(M)$ in the same way as above. Since $\tr^N_M(g)\in M_k^!(M)$ and the $g|_k\gamma$ vanish at every cusp of $\Gamma_0(N)$ that is not equivalent to $\infty$ under the action of $\Gamma_0(M)$, we see that $\tr^N_M(g)\in \mkh(M)$.

By an argument similar to the previous case, by considering $\tr^N_M(g)$ as a modular form in the larger space $M_k^!(N)$ we find that the principal part and constant term of $g$ determine the principal part and constant term of the Fourier expansion of $\tr^N_M(g)$. Because $\tr^N_M(g)\in \mkh(M)$ and $S_k(M)=0$, the form $\tr^N_M(g)$ is determined completely by its principal part and the constant term of its Fourier expansion. Thus, $\tr^N_M(g)$ is equal to the linear combination of basis elements from $\mkh(M)$ with principal part and constant term that equal the principal part and constant term of $g$.

\end{proof}

\subsection{Example: Traces from Level 4}
Let $k\leq 2$ with $ k\neq 0$ and let $f\in \mki(4)$. Note that $\Gamma_0(4)$ has three cusps which we take as $\infty,$ $0$, and $\frac{1}{2}$. Then
\begin{equation*}
    \tr^4_1(f)=\sum_{\gamma\in \Gamma_0(4)\backslash\slz}f|_k\gamma=f+f|_k\begin{psmallmatrix}
        0&-1\\
        1&0\\
    \end{psmallmatrix}+
    f|_k \begin{psmallmatrix}
        0&-1\\
        1&1\\
    \end{psmallmatrix}+
    f|_k\begin{psmallmatrix}
        0&-1\\
        1&2\\
    \end{psmallmatrix}+
    f|_k\begin{psmallmatrix}
        0&-1\\
        1&3\\
    \end{psmallmatrix}+
    f|_k\begin{psmallmatrix}
        1&0\\
        2&1\\
    \end{psmallmatrix}.
\end{equation*}
All of $f|_k\begin{psmallmatrix}
        0&-1\\
        1&0\\
    \end{psmallmatrix},
    f|_k \begin{psmallmatrix}
        0&-1\\
        1&1\\
    \end{psmallmatrix},
    f|_k\begin{psmallmatrix}
        0&-1\\
        1&2\\
    \end{psmallmatrix},$ and $
    f|_k\begin{psmallmatrix}
        0&-1\\
        1&3\\
    \end{psmallmatrix}$ can be interpreted as Fourier expansions of $f$ about the cusp 0. Because $f$ is holomorphic at 0, we know that these forms are all holomorphic at $\infty$. Similarly, since $f$ is holomorphic at $\frac{1}{2}$, we know that $f|_k\begin{psmallmatrix}
        1&0\\
        2&1\\
    \end{psmallmatrix}$ is holomorphic at $\infty$. Hence, as $\tr^4_1(f)\in M_k^!(1)$, and there are no holomorphic forms in $M_k^!(1)$ for these $k$, the principal part of $f$ completely determines $\tr^4_1(f)$.

    Now consider the trace from level 4 to level 2.
    Let $f\in\mki(4)$ and let $k\leq-2$. Then
    \begin{equation*}
    \tr^4_2(f)=\sum_{\gamma\in \Gamma_0(4)\backslash\Gamma_0(2)}f|_k\gamma=f+
    f|_k\begin{psmallmatrix}
        1&0\\
        2&1\\
    \end{psmallmatrix}.
\end{equation*}
We know that $\tr^4_2(f)\in M_k^!(2)$, but to calculate it explicitly, we show that $\tr^4_2(f)\in\mki(2)$.
The form $f|_k\begin{psmallmatrix}
        1&0\\
        2&1\\
    \end{psmallmatrix}$
is an expansion of $f$ about the cusp $\frac{1}{2}$. Hence, it is holomorphic at $\infty$ because $f$ is holomorphic at $\frac{1}{2}$. To show that it is also holomorphic at $0$, consider $f|_k\begin{psmallmatrix}
        1&0\\
        2&1\\
    \end{psmallmatrix}|_k\begin{psmallmatrix}
        0&-1\\
        1&0\\
    \end{psmallmatrix}=f|_k\begin{psmallmatrix}
        0&-1\\
        1&-2\\
    \end{psmallmatrix}$.
   This is an expansion of $f|_k\begin{psmallmatrix}
        1&0\\
        2&1\\
    \end{psmallmatrix}$ about 0 and also an expansion of $f$ about 0. Thus, since $f$ is holomorphic at 0, so is $f|_k\begin{psmallmatrix}
        1&0\\
        2&1\\
    \end{psmallmatrix}$. This fact, combined with the fact that $\infty$ and $\frac{1}{2}$ are equivalent under $\Gamma_0(2)$, means that $\tr^4_2(f)\in\mki(2)$. Since there are no holomorphic forms in $\mki(2)$ for $k<0$, the principal part of $f$ completely determines $\tr^4_2(f)$.

\section{Proof of Theorem 1.4}
 We now use the results from the previous section to prove Theorem \ref{MainTheorem}.
\begin{customthm}{1.4}\label{MainTheoremRestate}
Let $N\in\{2, 3, 4, 5, 6, 7, 8, 9, 10, 12, 13, 16, 18, 25\}$, so that $\Gamma_0(N)$ has genus zero. Let $k$ be an even integer and let $M\mid N$. Then the trace operator $\tr^N_M$ preserves duality between grids of basis elements for $\mki(N)$ and $\mkhd(N)$ for $k=0$. Additionally, if $N\in\{2, 3, 4\}$ then $\tr^N_M$ preserves duality for $k=-2$, and $\tr^2_1$ also preserves duality for $k=-4$. If $M=N,$ then $\tr^N_M$ is the identity map and preserves duality. In all other weights and levels $\tr^N_M$ does not preserve duality.
\end{customthm}

\begin{proof}
Let $N\in\{2, 3, 4, 5, 6, 7, 8, 9, 10, 12, 13, 16, 18, 25\}$, let $k\in 2\z$, and let $M\mid N$ with $M>0$. If $M=N$, then $\tr^N_M$ is the identity map, and it trivially preserves duality. Thus, suppose that $M\neq N$.  We note that by the formulas for $v(N)$ and $u(N)$ given in the appendix, if $M\mid N$, then $|v(N)|\geq |v(M)|$ and $|u(N)|\geq |u(M)|$.

We break the proof into three cases. The first is when duality is not preserved due to obstructions that arise from applying the trace to the $f\in\mki(N)$, the second is when duality is not preserved due to obstructions that arise from applying the trace to the $g\in\mkhd(N)$, and the third is when no obstructions occur from either space.

   In the first case, suppose that $v_k(N)<v_k(M)$ and $v_k(N)<0$. In particular, this occurs when $k<-4$ for $N=2$, when $k<-2$ for $N=3,4$, and when $k<0$ for all other $N$. Let $f^{(N)}_{k,m}(z)$ be a canonical basis element of $\mki(N)$. By Proposition \ref{Prop:LvN Trace}, we know $\tr^N_M(f^{(N)}_{k,m}(z))$ is determined by the principal part of $f^{(N)}_{k,m}(z)$. The Fourier expansion of the form $f^{(N)}_{k,m}(z)$ may be written
    \begin{equation}
        f^{(N)}_{k,m}(z)=q^{-m}+a_k^{(N)}(m,v(N)+1)q^{v(N)+1}+\cdots+a_k^{(N)}(m,v(M))q^{v(M)}+O\left(q^{v(M)+1}\right).
    \end{equation}
    For each $n\leq v(M)$ there is a canonical basis element of $\mki(M)$ with Fourier expansion of the form
    \begin{equation}
        f^{(M)}_{k,-n}(z)=q^n+O\left(q^{v(M)+1}\right).
    \end{equation}
    Since the principal part of $\tr^N_M(f_{k,m}^{(N)})$ equals the principal part of $f^{(N)}_{k,m}$, we may write the trace in terms of canonical basis elements as
    \begin{equation}\label{eq:Trace f}
        \tr^N_M(f^{(N)}_{k,m}(z))=f^{(M)}_{k,m}(z)+a_k^{(N)}(m,v(N)+1)f^{(M)}_{k,-(v(N)+1)}(z)+\cdots+a_k^{(N)}(m,v(M))f^{(M)}_{k,-v(M)}(z).
    \end{equation}
Applying the trace to one basis element preserves the terms of its Fourier expansion with exponents less than $v(M)+1$, and applying the trace to the grid of all basis elements preserves all columns of Fourier coefficients with exponents less than $v(M)+1$. Thus, if the trace preserves duality, then there must be a modular form $g(z)\in\mkhd(M)$ with Fourier expansion
\begin{equation}
    g(z)=q^{-v(N)-1}-\sum_{n=m}^\infty a_k^{(N)}(n,v(N)+1)q^{n}.
\end{equation}
Applying duality in level $N$ and comparing Fourier expansions, we find that if $g(z)$ exists, it must equal the basis element $g^{(N)}_{2-k,v(N)+1}(z)$ in the space $\mkhd(N)$, and must also be a form of level $M$.

If $g(z)$ is of level $M$ and weight $2-k$, its order of vanishing at $\infty$, which is $-v_k(N)-1$, is bounded above by $u_{2-k}(M)$. Formulas for $v_k(N)$ and $u_{2-k}(M)$ in the appendix show that in all cases where $-v_k(N)>-v_k(M)$ and $v(N)<0$, we have $-v_k(N)-1>u_{2-k}(M)$, a contradiction. Therefore, the trace cannot preserve duality in this case.

   In the second case, suppose that $u_{2-k}{(N)}<u_{2-k}{(M)}$ and $u_{2-k}{(N)}<1$. This occurs when $2-k<2$, so $k>0$. Let $g^{(N)}_{2-k,n}(z)\in\mkhd(N)$. Then, by Proposition \ref{Prop:LvN Trace}, the principal part and constant term of $g^{(N)}_{2-k,n}(z)$ completely determine $\tr^N_M(g^{(N)}_{2-k,n}(z))$.
    Thus, if
    \begin{equation}
        g^{(N)}_{2-k,n}(z)=q^{-n}+b_{2-k}^{(N)}(n,u(N)+1)q^{u(N)+1}+\cdots+b_{2-k}^{(N)}(n,u(M))q^{u(M)}+O\left(q^{u(M)+1}\right),
    \end{equation}
    then by using our basis for the space $\mkhd(M)$, we find that
    \begin{equation}\label{eq:Trace g}
        \tr^N_M(g^{(N)}_{2-k,n}(z))=g^{(M)}_{2-k,n}(z)+b_{2-k}^{(N)}(n,u(N)+1)g^{(M)}_{2-k,u(N)+1}(z)+\cdots+b_{2-k}^{(N)}(n,u(M))g^{(M)}_{2-k,u(M)}(z).
    \end{equation}
    Here, applying the trace to one basis element preserves the terms of its Fourier expansion with exponents less than $u(M)+1$. Taking the trace of all basis elements in the grid preserves all columns of Fourier coefficients with exponents less than $u(M)+1$. Similar to the previous case, this means that there must be a modular form $f(z)\in\mki(M)$ with Fourier expansion
    \begin{equation}
        f(z)=q^{-u(N)-1}-\sum_{m=n}^\infty b_{2-k}^{(N)}(m,u(N)+1)q^m.
    \end{equation}
It follows that, if $f(z)$ exists, it must equal the basis element $f^{(N)}_{k,u(N)+1}(z)\in\mki(N)$, and be a modular form of level $M$.
If $f(z)$ is of level $M$ and weight $k$, its order of vanishing at infinity, $-u_{2-k}(N)-1$, is bounded above by $v_k(M)$. Formulas for $u_{2-k}(N)$ and $v_k(M)$ in the appendix show that when $u_{2-k}{(N)}<u_{2-k}{(M)}$ and $u_{2-k}{(N)}<1$, we have  $-u_{2-k}(N)-1>v_k(M)$, a contradiction. Thus, in this case the trace cannot preserve duality.

    In the remaining spaces, the trace operator does preserve duality. We note that this is a special case of \cite[Theorem~7.4]{griffin}. Suppose that $N,M$ are as before, but that $v(N)=v(M)$ with $v(N)<0$, and $u(N)=u(M)$ with $u(N)<1$. This occurs for any $N$ when $k=0$, when $N=2,3,$ or $4$ and $k=-2$, and when $N=2$ and $k=-4$. Let $f^{(N)}_{k,m}\in\mki(N)$ and $g^{(N)}_{2-k,n}\in\mkhd(N)$ be canonical basis elements. Again, by Proposition \ref{Prop:LvN Trace}, we find that $\tr^N_M(f^{(N)}_{k,m})$ is determined by the principal part of $f^{(N)}_{k,m}$, but now $f^{(N)}_{k,m}$ and $f^{(M)}_{k,m}$ have identical principal part. Similarly, $\tr^N_M(g^{(N)}_{2-k,n})$ is determined by the principal part and constant term of $g^{(N)}_{2-k,n}$, and $g^{(N)}_{2-k,n}$ and $g^{(M)}_{2-k,n}$ have identical principal part. Hence,
    \begin{equation}
       \tr^N_M(f^{(N)}_{k,m})=f^{(M)}_{k,m} \quad\textrm{and}\quad\tr^N_M(g^{(N)}_{2-k,n})=g^{(M)}_{2-k,n}.
    \end{equation}
    This means that the traces of the canonical basis elements in level $N$ are the canonical basis elements in level $M$ that are already form a modular grid. Thus, the trace preserves duality.
\end{proof}

\section{Generating Functions}
Generating functions for modular grids and the action of linear operators on modular grids are described in \cite[Section~7]{griffin}. These generating functions can describe the obstructions that occur when duality between grids is not preserved. We denote the generating function for the modular grid for $\mki(N)$ and $\mkhd(N)$ by $\mathcal{H}_{k,N}(z,\tau)$. Let $p=e^{2\pi i \tau}$; then by Theorem $\ref{GenFunTheorem}$ for appropriate regions of $\mathcal{H}\times\mathcal{H}$ we have that
\[\mathcal{H}_{k,N}(z,\tau)=\sum_{m=-v_k(N)}^\infty f_{k,m}^{(N)}(z)p^m \quad\textrm{and}\quad \mathcal{H}_{k,N}(z,\tau)=-\sum_{n=-u_{2-k}(N)}^\infty g_{2-k,n}^{(N)}(\tau)q^n.\]
Additionally, let $\tr^N_M(\mathcal{H}_{k,N}(z,\tau))_{k,z}$ denote the trace applied to each canonical basis element in the generating function in weight $k$ and variable $z$, and recall from \eqref{eq:basisFourier} that $a_k^{(N)}(m,n)$ and $b_k^{(N)}(m,n)$ are the Fourier coefficients of basis elements from $\mki(N)$ and $\mkh(N)$ respectively.

As an example, the following forms appear in the modular grids for $M_{-6}^{(\infty)}(2)$ and $\widehat{M}_8^{(\infty)}(2)$:
    \begin{align*}
        f_{-6,2}^{(2)}(z)&= q^{-2}+8q^{-1}-224+2144q+\cdots, \\
        f_{-6,3}^{(2)}(z)&= q^{-3}-12q^{-1}+4096-98226q+\cdots, \\
        f_{-6,4}^{(2)}(z)&= q^{-4}-64q^{-1}-31200+1817856q+\cdots, \\
        &\\
        g_{8,-1}^{(2)}(z)&=q-8q^2+12q^3+64q^4+\cdots, \\
        g_{8,0}^{(2)}(z)&=1+224q^2-4096q^3+31200q^4+\cdots,\\
        g_{8,1}^{(2)}(z)&=q^{-1}-2144q^2+98226q^3-181756q^4+\cdots.
    \end{align*}
     Upon calculating their traces, we find that
    \begin{align*}
        \tr^2_1(f_{-6,2}^{(2)}(z))&=f_{-6,2}^{(1)}(z)+8f_{-6,1}^{(1)}(z)= q^{-2}+8q^{-1}-65760-87553952q+\cdots,\\
        \tr^2_1(f_{-6,3}^{(2)}(z))&=f_{-6,3}^{(1)}(z)-12f_{-6,1}^{(1)}(z)= q^{-3}-12q^{-1}-1044480-22875832242q+\cdots,\\
        \tr^2_1(f_{-6,4}^{(2)}(z))&=f_{-6,4}^{(1)}(z)-64f_{-6,1}^{(1)}(z)= q^{-4}-64q^{-1}-7895520-1969010000640q+\cdots,\\
        &\\
        \tr^2_1(g_{8,-1}^{(2)}(z))&=0,\\
        \tr^2_1(g_{8,0}^{(2)}(z))&=f_{8,0}^{(1)}(z)=1+480q+61920q^2+1050240q^3+\cdots,\\
        \tr^2_1(g_{8,1}^{(2)}(z))&=f_{8,1}^{(1)}(z)=q^{-1}+28240q+87326720q^2+22876173090q^3+\cdots.
    \end{align*}
    Inspecting the Fourier coefficients immediately shows that in this case, duality is not preserved by the trace operator.

    Explicitly computing the action of the trace on the generating function in weight $8$ and variable $\tau$ gives
    \begin{align*}
        \tr^2_1(\mathcal{H}_{-6,2}(z,\tau))_{8,\tau}&=\sum_{n=-1}^\infty\tr^2_1(g^{(2)}_{8,n})(\tau)q^n\\
        &=\sum_{n=0}^\infty g^{(1)}_{8,n}(\tau)q^n\\
          &=\mathcal{H}_{-6,1}(z,\tau).
    \end{align*}
    In weight $-6$ and variable $z$, however, the action of the trace is
    \begin{align*}
         \tr^2_1(\mathcal{H}_{-6,2}(z,\tau))_{-6,z}&=\tr^2_1\left(\sum_{m=2}^\infty f_{-6,m}^{(2)}(z)p^m\right)\\
         &=\sum_{m=2}^\infty\tr^2_1\left( f_{-6,m}^{(2)}(z)\right)p^m\\
         &=\sum_{m=2}^\infty\left( f_{-6,m}^{(1)}(z)+a_{-6}^{(2)}(m,1)f_{-6,1}^{(1)}(z)\right)p^m\\
         &=\sum_{m=2}^\infty f_{-6,m}^{(1)}(z)p^m-\sum_{m=2}^\infty b_{8}^{(2)}(1,m)f_{-6,1}^{(1)}(z)p^m\\
         &=\sum_{m=1}^\infty f_{-6,m}^{(1)}(z)p^m-\sum_{m=1}^\infty b_{8}^{(2)}(1,m)f_{-6,1}^{(1)}(z)p^m\\
         &=\mathcal{H}_{-6,1}(z,\tau)-f^{(1)}_{-6,1}(z)g^{(2)}_{8,-1}(\tau).
    \end{align*}

More generally, when $v_k(N)< v_k(M)$ and $v_k(N)<0$, by equation \eqref{eq:Trace f} we have
\begin{align*}
    \tr^N_M(\mathcal{H}_{k,N}(z,\tau))_{k,z}=&\sum_{m=-v_k(N)}^\infty \tr^N_M(f_{k,m}^{(N)}(z))p^m\\
    =&\sum_{m=-v_k(N)}^\infty\biggl(f^{(M)}_{k,m}(z)+a_k^{(N)}(m,v_k(N)+1)f^{(M)}_{k,-(v_k(N)+1)}(z)+\cdots\\
    &+a_k^{(N)}(m,v_k(M))f^{(M)}_{k,-v_k(M)}(z)\biggr)p^m\\
    =&\sum_{m=-v_k(N)}^\infty f_{k,m}^{(M)}(z)p^m-\sum_{m=-v_k(N)}^\infty b_{2-k}^{(N)}(v_k(N)+1,m)f_{k,-v_k(N)-1}^{(M)}(z)p^m-\cdots\\
    &-\sum_{m=-v_k(N)}^\infty b_{2-k}^{(N)}(v_k(M),m)f^{(M)}_{k,-v_k(M)}(z)p^m\\
    =&\mathcal{H}_{k,M}(z,\tau)-\sum_{x=1}^{v_k(M)-v_k(N)}f^{(M)}_{k,-(v_k(N)+x)}(z)g^{(N)}_{2-k,(v_k(N)+x)}(\tau).
\end{align*}
Here the last equality comes because $b_{2-k}(n,m)=0$ when $m<-n$. If $v(N)=v(M)$ we have that
\begin{align*}
    \tr^N_M(\mathcal{H}_{k,N}(z,\tau))_{k,z}=&\sum_{m=-v_k(N)}^\infty \tr^N_M(f_{k,m}^{(N)}(z))p^m\\
    &=\sum_{m=-v_k(N)}^\infty f^{(M)}_{k,m}(z)p^m\\
    &=\mathcal{H}_{k,M}(z,\tau).
\end{align*}
Hence, when $v_k(N)\leq v_k(M)$ and $v_k(N)<0$, taking traces of the grid in weight $k$ gives
 \begin{equation} \label{eq:gen1}
        \tr^N_M(\mathcal{H}_{k,N}(z,\tau))_{k,z}=\mathcal{H}_{k,M}(z,\tau)-\sum_{x=1}^{v(M)-v(N)}f^{(M)}_{k,-(v(N)+x)}(z)g^{(N)}_{2-k,(v(N)+x)}(\tau).
    \end{equation}
Similarly, when $u_{2-k}(N)\leq u_{2-k}(M)$ and $u_{2-k}(N)<1$, by equation \eqref{eq:Trace g} taking traces of the grid in weight $2-k$ gives
 \begin{equation} \label{eq:gen2}
        \tr^N_M(\mathcal{H}_{k,N}(z,\tau))_{2-k,\tau}=\mathcal{H}_{k,M}(z,\tau)-\sum_{x=1}^{u(M)-u(N)}f^{(N)}_{k,(u(N)+x)}(z)g^{(M)}_{2-k,-(u(N)+x)}(\tau).
    \end{equation}
When the sums in $\eqref{eq:gen1}$ and $\eqref{eq:gen2}$ are both empty, it is clear that the trace operator preserves duality.

\appendix
\section*{Appendix: Building Blocks for Canonical Bases}
Here we list the modular forms used to construct the basis for each space where $\Gamma_0(N)$ has genus zero. The modular forms for level 1 are from \cite{duke}, those of level 2 (and the method for finding those of level 3) are from \cite{garthwaite}, the method of finding the elements for the remaining prime levels can be found in \cite{thorntonp}, levels 6, 10, 12, and 18 can be found in \cite{Iba}, and levels 8, 9, 16, and 25 can be found in \cite{thornton}. We correct any typos where applicable. For information about which spaces of modular forms are spanned by eta quotients see \cite{ROUSE}. Our calculations were performed using Sage \cite{sagemath}.

\subsection*{Level 1}
In level 1, we write $k=12\ell+k',$ with $ k'\in\{0,4,6,8,10,14\}$ and $v=\ell=u$. There is one cusp in level 1, which is taken to be $\infty$. The first basis element of $\mki(1)$ is written as
\begin{equation}
    f_{k,-v}^{(1)}(z)=\Delta(z)^\ell E_{k'}(z),
\end{equation}
where $\Delta$ is the normalized cusp form of weight 12 and level 1, $E_0=1$, and otherwise $E_{k'}(z)$ is the Eisenstein series of weight $k'$. We choose the Hauptmodul
\begin{equation*}
   j(z)=q^{-1}+744+196884q+21493760q^2+864299970q^3+\cdots,
\end{equation*}
the $j$-function. See \cite{duke} for further details.

\subsection*{Level 2}
In level 2, we write $k=4\ell+k',$ with $k'\in\{0,2\}$ so $v=\ell$ and $u=v-1$. There are two cusps in level 2, which are taken to be $\infty, 0$. The first basis element of $\mki(2)$ is
\begin{equation}
    f^{(2)}_{k,-v}(z)=F_4(z)^\ell F_{k'}(z),
\end{equation}
and is computed using
\begin{align}
    F_0(z)&=1,\\
    F_2(z)&=f^{(2)}_{2,0}(z)=2E_2(2z)-E_2(z)=1+24q+24q^2+96q^3+24q^4+\cdots,\\
    F_4(z)&=f^{(2)}_{4,-1}(z)=\frac{E_4(z)-E_4(2z)}{240}=q+8q^2+28q^3+64q^4+\cdots.
\end{align}
A Hauptmodul is given by
\begin{equation}
    \psi^{(2)}(z)=\frac{\eta(z)^{24}}{\eta(2z)^{24}}=q^{-1}-24+276q-2048q^2+\cdots.
\end{equation}
Its value at the cusp $0$ is 0.

\subsection*{Level 3}
In level 3, we write $k=6\ell+k',$ with $k'\in\{0,2,4\}$ so $v=2\ell+\lfloor\frac{k'}{3}\rfloor$ and $u=v-1$. There are two cusps in level 3, which are taken to be $\infty, 0$. The first basis element of $\mki(3)$ is
\begin{equation}
    f^{(3)}_{k,-v}(z)=F_6(z)^\ell F_{k'}(z),
\end{equation}
and is computed using
\begin{align}
    F_0(z)&=1,\\
    F_2(z)&=f^{(3)}_{2,0}(z)=\frac{3E_2(3z)-E_2(z)}{2}=1+12q+36q^2+12q^3+84q^4+\cdots,\\
   F_4(z)&=\displaystyle f^{(3)}_{4,-1}(z)=\frac{E_4(z)-\left(f^{(3)}_{2,0}(z)\right)^2}{216}=q+9q^2+27q^3+73q^4+126q^5+\cdots,\\
F_6(z)&=f^{(3)}_{6,-2}(z)=\frac{\eta(3z)^{18}}{\eta(z)^6}q^2+6q^3+27q^4+80q^5+207q^6+\cdots.
\end{align}
A Hauptmodul is given by
\begin{equation}
\psi^{(3)}(z)=\frac{\eta(z)^{12}}{\eta(3z)^{12}}=q^{-1}-12+54q-76q^2+\cdots.
\end{equation}
Its value at the cusp $0$ is 0.
\subsection*{Level 4}
In level 4, $v=\frac{k}{2}$ and $u=v-2$. There are three cusps in level 4, which are taken to be $0,\frac{1}{2},\infty$. The first basis element of $\mki(4)$ is
\begin{equation}
    f^{(4)}_{k,-v}(z)=f^{(4)}_{2,-1}(z)^{k/2},
\end{equation}
where
\begin{align}
  f^{(4)}_{2,-1}(z)=3E_2(2z)-E_2(z)-2E_2(4z)&=q+4q^3+6q^5+8q^7+13q^9+\cdots.
\end{align}
A Hauptmodul is given by
\begin{equation}
\psi^{(4)}(z)=\frac{\eta(z)^8}{\eta(4z)^8}=q^{-1}-8+20q-62q^2+\cdots.
\end{equation}
Its values at the cusps $0$ and $\frac{1}{2}$ are 0 and $-16$.
\subsection*{Level 5}
In level 5, we write $k=4\ell+k',$ with $k'\in\{0,2\}$ so $v=2\ell$ and $u=v-1$. There are two cusps in level 5, which are taken to be $\infty, 0$. The first basis element of $\mki(5)$ is
\begin{equation}
    f^{(5)}_{k,-v}(z)=F_4(z)^\ell F_{k'}(z),
\end{equation}
and is computed using
\begin{align}
F_0(z)&=1,\\
    F_2(z)&=f^{(5)}_{2,0}(z)=\frac{5E_2(5z)-E_2(z)}{4}=1+6q+18q^2+24q^3+42q^4+\cdots,\\
F_4(z)&=f^{(5)}_{4,-2}(z)=\frac{\eta(5z)^{10}}{\eta(z)^2}=q^2+2q^3+5q^4+10q^5+20q^6+\cdots.
\end{align}
A Hauptmodul is given by
\begin{equation}
\psi^{(5)}(z)=\frac{\eta(z)^6}{\eta(5z)^6}=q^{-1}-6+9q+10q^2-30q^3+\cdots.
\end{equation}
Its value at the cusp $0$ is 0.
\subsection*{Level 6}
In level 6, $v=k$ and $u=v-3$. There are four cusps in level 6, which are taken to be $0,\frac{1}{3},\frac{1}{2},\infty$. The first basis element of $\mki(6)$ is
\begin{equation}
    f^{(6)}_{k,-v}(z)=f^{(6)}_{2,-2}(z)^{k/2},
\end{equation}
where
\begin{align}
  f^{(6)}_{2,-2}(z)=\frac{\eta(z)^2\eta(6z)^{12}}{\eta(2z)^4\eta(3z)^6}=q^2-2q^3+3q^4-q^6+7q^8+\cdots.
\end{align}
A Hauptmodul is given by
\begin{equation}
\psi^{(6)}(z)=\frac{\eta(2z)^8\eta(3z)^4}{\eta(z)^4\eta(6z)^8}=q^{-1}+4+6q+4q^2-3q^3+\cdots.
\end{equation}
Its values at the cusps away from $\infty$ are the roots of the polynomial $x^3-10x+9x$.
\subsection*{Level 7}
In level 7, we write $k=6\ell+k',$ with $k'\in\{0,2,4\}$ so $v=4\ell+2\lfloor\frac{k'}{3}\rfloor$ and $u=v-1$. There are two cusps in level 7, which are taken to be $\infty, 0$. The first basis element of $\mki(7)$ is
\begin{equation}
    f^{(7)}_{k,-v}(z)=F_6(z)^\ell F_{k'}(z),
\end{equation}
and is computed using
\begin{align}
F_0(z)&=1,\\
   F_2(z)&=f^{(7)}_{2,0}(z)=\frac{7E_2(7z)-E_2(z)}{6}=1+4q+12q^2+16q^3+28q^4+\cdots,\\
    F_4(z)&=f^{(7)}_{4,-2}(z)=q^2+3q^3+8q^4+11q^5+\cdots,\\
    F_6(z)&=f^{(7)}_{6,-4}(z)=\frac{\eta(7z)^{14}}{\eta(z)^2}=q^4+2q^5+5q^6+10q^7+20q^8+\cdots.
\end{align}
Level 7 is not spanned by eta-quotients, so we use Sage to calculate the Fourier expansion of $f^{(7)}_{4,-2}$.
A Hauptmodul is given by
\begin{equation}
\psi^{(7)}(z)=\frac{\eta(z)^4}{\eta(7z)^4}=q^{-1}-4+2q+8q^2-5q^3+\cdots.
\end{equation}
Its value at the cusp $0$ is 0.
\subsection*{Level 8}
In level 8, $v=k$ and $u=v-3$. There are four cusps in level 8, which are taken to be $0,\frac{1}{4},\frac{1}{2},\infty$. The first basis element of $\mki(8)$ is
\begin{equation}
    f^{(8)}_{k,-v}=f^{(8)}_{2,-2}(z)^{k/2},
\end{equation}
where
\begin{align}
 f^{(8)}_{2,-2}(z)=\frac{\eta(8z)^{8}}{\eta(4z)^4}&=q^2+4q^6+6q^{10}+8q^{14}+13q^{18}+\cdots.
\end{align}
A Hauptmodul is given by
\begin{equation}
\psi^{(8)}(z)=\frac{\eta(z)^4\eta(4z)^2}{\eta(2z)^2\eta(8z)^4}=q^{-1}-4+4q+2q^3-8q^5+\cdots.
\end{equation}
Its values at the cusps away from $\infty$ are the roots of the polynomial $x^3+12x^2+32x$.
\subsection*{Level 9}
In level 9, $v=k$ and $u=v-3$. There are four cusps in level 9, which are taken to be $0,\pm\frac{1}{3},\infty$. The first basis element of $\mki(9)$ is
\begin{equation}
    f^{(9)}_{k,-v}(z)=f^{(9)}_{2,-2}(z)^{k/2},
\end{equation}
where
\begin{align}
f^{(9)}_{2,-2}(z)=\frac{\eta(9z)^{6}}{\eta(3z)^2}&=q^2+2q^5+5q^8+4q^{11}+8q^{14}+\cdots.
\end{align}
A Hauptmodul is given by
\begin{equation}
\psi^{(9)}(z)=\frac{\eta(z)^4\eta(4z)^2}{\eta(2z)^2\eta(8z)^4}=q^{-1}-4+4q+2q^3-8q^5+\cdots.
\end{equation}
Its values at the cusps away from $\infty$ are the roots of the polynomial $x^3+9x^2+27x$.
\subsection*{Level 10}
In level 10, we write $k=4\ell+k',$ with $k'\in\{0,2\}$ so $v=6\ell+k'$ and $u=v-3$. There are four cusps in level 10, which are taken to be $0,\frac{1}{5},\frac{1}{2},\infty$. The first basis element of $\mki(10)$ is
\begin{equation}
    f^{(10)}_{k,-v}(z)=F_4(z)^\ell F_{k'}(z),
\end{equation}
and is computed using
\begin{align}
F_0(z)&=1,\\
    F_2(z)&=f^{(10)}_{2,-2}(z)=q^2+3q^4-4q^5+4q^6+7q^8+\cdots,\\
F_4(z)&=f^{(10)}_{4,-6}(z)=q^6-2q^7+3q^8-6q^9+11q^{10}+\cdots.
\end{align}
Level 10 is not spanned by eta-quotients, so we use Sage to calculate the Fourier expansions of $f^{(10)}_{2,-2}$ and $f^{(10)}_{4,-6}$. Note that these are not the forms in \cite{Iba} for level 10, weights 2 and 4, which are incorrect.

A Hauptmodul is given by
\begin{equation}
\psi^{(10)}(z)=\frac{\eta(2z)\eta(5z)^5}{\eta(z)\eta(10z)^5}=q^{-1}+1+q+2q^2+2q^3+\cdots.
\end{equation}
Its values at the cusps away from $\infty$ are the roots of the polynomial $x^3-3x^2-4x$.
\subsection*{Level 12}
In level 12, $v=2k$ and $u=v-5$. There are six cusps in level 12, which are taken to be $0,\frac{1}{6},\frac{1}{4},\frac{1}{3},\frac{1}{2},\infty$. The first basis element of $\mki(12)$ is
\begin{equation}
      f^{(12)}_{k,-v}(z)=f^{(12)}_{2,-4}(z)^{k/2},
\end{equation}
where
\begin{align}
    f^{(12)}_{2,-4}(z)=&\frac{\eta(z)^{10}\eta(4z)\eta(6z)^9}{27\eta(2z)^7\eta(3z)^6\eta(12z)^3}+\frac{11\eta(z)^{7}\eta(4z)^4\eta(6z)^9}{72\eta(2z)^7\eta(3z)^5\eta(12z)^4}-\frac{\eta(z)^{4}\eta(4z)^7\eta(6z)^9}{12\eta(2z)^7\eta(3z)^4\eta(12z)^5}\\
    &+\frac{\eta(z)\eta(4z)^{10}\eta(6z)^9}{54\eta(2z)^7\eta(3z)^3\eta(12z)^6}-\frac{\eta(z)^{9}\eta(4z)^3\eta(6z)^2}{8\eta(2z)^6\eta(3z)^3\eta(12z)}.
\end{align}
This form has Fourier expansion given by
\begin{equation}
    f^{(12)}_{2,-4}(z)=q^4-2q^6+3q^8-q^{12}+7q^{16}+\cdots.
\end{equation}
A Hauptmodul is given by
\begin{equation}
\psi^{(12)}(z)=\frac{\eta(4z)^4\eta(6z)^2}{\eta(2z)^2\eta(12z)^4}=q^{-1}+2q+q^3-2q^7+\cdots.
\end{equation}
Its values at the cusps away from $\infty$ are the roots of the polynomial $x^5-10x^3+9x$.
\subsection*{Level 13} In level 13, we write $k=12\ell+k'$ with $ k'\in\{0,2,4,6,8,10\}$.  We find that $v=14\ell+k'$, unless $k'=2$, in which case $v=14\ell$, and $u=v-1$. There are two cusps in level 13, which are taken to be $\infty, 0$. The first basis element of $\mki(13)$ is
\begin{equation}
    f^{(13)}_{k,-v}(z)=F_{12}(z)^\ell F_{k'}(z),
\end{equation}
and is computed using
\begin{align*}
F_0(z)&=1,\\
    F_2(z)&=f_{2,0}^{(13)}(z)=\frac{13E_2(13z)-E_2(z)}{12}=1+2q+6q^2+8q^3+14q^4+\cdots,\\
    F_4(z)&=f_{4,-4}^{(13)}(z)=q^4+q^5+q^6-q^7-3q^9+\cdots,\\
    F_6(z)&=f_{6,-6}^{(13)}(z)=q^6+q^7+q^8+3q^9+2q^{11}+\cdots,\\
    F_8(z)&=f_{8,-8}^{(13)}(z)=f_{4,-4}^{(13)}(z)^2=q^8+2q^9+3q^{10}-q^{12}-8q^{13}+\cdots,\\
   F_{10}(z)&= f_{10,-10}^{(13)}(z)= f_{4,-4}^{(13)}(z)f_{6,-6}^{(13)}(z)=q^{10}+2q^{11}+3q^{12}+4q^{13}+3q^{14}+\cdots,\\
    F_{12}(z)&=f^{(13)}_{12,-14}(z)=\frac{\eta(13z)^{26}}{\eta(z)^{2}}=q^{14}+2q^{15}+5q^{16}+10q^{17}+20q^{18}+\cdots.
\end{align*}
Level 13 is not spanned by eta-quotients, so we use Sage to calculate the Fourier expansions of $f^{(13)}_{4,-4}$ and $f^{(13)}_{6,-6}$. A Hauptmodul is given by
\begin{equation}
   \psi^{(13)}(z)=\frac{\eta(z)^2}{\eta(13z)^2}=q^{-1}-2-q+2q^2+q^3+\cdots.
\end{equation}
Its value at the cusp $0$ is 0.

\subsection*{Level 16}In level 16, $v=2k$ and $u=v-5$. There are six cusps in level 16, which are taken to be $0,\frac{1}{8},\pm\frac{1}{4},\frac{1}{2},\infty$. The first basis element of $\mki(16)$ is
\begin{equation}
    f^{(16)}_{k,-v}(z)=f^{(16)}_{2,-4}(z)^{k/2},
\end{equation}
where
\begin{equation}
    f^{(16)}_{2,-4}(z)=\frac{\eta(16z)^{8}}{\eta(8z)^{4}}=q^4+4q^{12}+6q^{20}+8q^{28}+13q^{36}+\cdots.
\end{equation}
A Hauptmodul is given by
\begin{equation}
    \psi^{(16)}(z)=\frac{\eta(z)^2\eta(8z)}{\eta(2z)\eta(16z)^2}=q^{-1}-2+2q^3-q^7+\cdots.
\end{equation}
Its values at the cusps away from $\infty$ are the roots of the polynomial $x^5+10x^4+40x^3+80x^2+64x$.
\subsection*{Level 18}In level 18, $v=3k$ and $u=v-7$. There are eight cusps in level 18, which are taken to be $0,\frac{1}{9},\pm\frac{1}{6},\pm\frac{1}{3},\frac{1}{2},\infty$.
The first basis element of $\mki(18)$ is
\begin{equation}
    f^{(18)}_{k,-v}(z)=f^{(18)}_{2,-6}(z)^{k/2},
\end{equation}
where
\begin{align}
    f^{(18)}_{2,-6}(z)=&\frac{25\eta(z)^8\eta(6z)^2\eta(9z)^4}{216\eta(2z)^4\eta(3z)^4\eta(18z)^2}-\frac{11\eta(z)^3\eta(6z)^8\eta(9z)^7}{144\eta(2z)^3\eta(3z)^6\eta(18z)^5}-\frac{121\eta(z)^6\eta(6z)^7\eta(9z)}{972\eta(2z)^3\eta(3z)^5\eta(18z)^2}\\
    &-\frac{41\eta(z)^6\eta(6z)^2\eta(9z)^6}{144\eta(2z)^3\eta(3z)^4\eta(18z)^3}+\frac{67\eta(z)^4\eta(6z)^7\eta(9z)^3}{144\eta(2z)^2\eta(3z)^5\eta(18z)^3}\\
    &+\frac{\eta(2z)^9\eta(3z)^8\eta(18z)}{972\eta(z)^6\eta(6z)^6\eta(9z)^2}-\frac{125\eta(z)\eta(2z)^4\eta(9z)^2}{1296\eta(3z)\eta(6z)\eta(18z)}.
\end{align}
This form has Fourier expansion given by
\begin{equation}
    f^{(18)}_{2,-6}(z)=q^6-2q^9+3q^{12}-q^{18}+7q^{24}+\cdots.
\end{equation}
A Hauptmodul is given by
\begin{equation}
    \psi^{(18)}(z)=\frac{\eta(6z)\eta(9z)^3}{\eta(3z)\eta(18z)^3}=q^{-1}+q^2+q^5-q^8+\cdots.
\end{equation}
Its values at the cusps away from $\infty$ are the roots of the polynomial $x^7-7x^4-8x$.

\subsection*{Level 25}In level 25, we write $k=4\ell+k',$ with $k'\in\{0,2\}$ so $v=10\ell+2k'$ and $u=v-5$. There are six cusps in level 25, which are taken to be $0,\pm\frac{1}{5},\pm\frac{2}{5},\infty$. The first basis element of $\mki(25)$ is
\begin{equation}
    f^{(25)}_{k,-v}(z)=F_4(z)^\ell F_{k'}(z),
\end{equation}
and is computed using
\begin{align}
F_0(z)&=1,\\
    F_2(z)&=f^{(25)}_{2,-4}(z)=q^{4}+q^{6}+2q^{9}+3q^{14}+2q^{16}+\cdots,\\
    F_4(z)&=f^{(25)}_{4,-10}(z)=\frac{\eta(25z)^{10}}{\eta(5z)^{2}}=q^{10}+2q^{15}+5q^{20}+10q^{25}+20q^{30}+\cdots.
\end{align}
Level 25 is not spanned by eta-quotients, so we use Sage to calculate the Fourier expansion of $f^{(25)}_{2,-4}$. A Hauptmodul is given by
\begin{equation}
    \psi^{(25)}(z)=\frac{\eta(z)}{\eta(25z)}=q^{-1}-1-q+q^4+q^6+\cdots.
\end{equation}
Its values at the cusps away from $\infty$ are the roots of the polynomial $x^5+5x^4+15x^3+25x^2+25x$.

\bibliographystyle{plain}

\end{document}